\documentclass[sn-mathphys]{sn-jnl}

\jyear{2021}%

\theoremstyle{thmstyleone}%
\newtheorem{theorem}{Theorem}
\theoremstyle{thmstyletwo}%
\newtheorem{example}{Example}%

\theoremstyle{thmstylethree}%

\begin{document}

\title[Na\"{i}ve estimator of Poisson regression with errors]{The na\"{i}ve estimator of a Poisson regression model with measurement errors}


\author[1]{\fnm{Kentarou} \sur{Wada}}\email{wadaken4269@gmail.com}

\author*[1]{\fnm{Takeshi} \sur{Kurosawa}}\email{tkuro@rs.tus.ac.jp}


\affil[1]{\orgdiv{Department of Applied Mathematics}, \orgname{Tokyo University of Science}, \orgaddress{\street{Kagurazaka 1-3}, \city{Shinjyuku-ku}, \postcode{1628601}, \state{Tokyo}, \country{Japan}}}




\abstract{We generalize the na\"{i}ve estimator of a Poisson regression model with measurement errors as discussed in Kukush et al. \cite{three}. The explanatory variable is not always normally distributed as they assume. In this study, we assume that the explanatory variable and measurement error are not limited to a normal distribution. We clarify the requirements for the existence of the na\"{i}ve estimator and derive its asymptotic bias and asymptotic mean squared error (MSE). In addition, we propose a consistent estimator of the true parameter by correcting the bias of the na\"{i}ve estimator. As illustrative examples, we present simulation studies that compare the performance of the na\"{i}ve estimator and new estimator for a Gamma explanatory variable with a normal error or a Gamma error.}




\maketitle

\section{Introduction}\label{S1}
We often cannot measure explanatory variables correctly in regression models because an observation may not be performed properly. The estimation result may be distorted when we estimate the model from data with measurement errors. We call models with measurement errors in an explanatory variable Error in Variable (EIV) models. In addition, actual phenomena often cannot be explained adequately by a simple linear structure, and the estimation of non-linear models, especially generalized linear models, from data with errors is a significant problem. Various studies have focused on non-linear EIV models (see, for example, \cite{box,gea}). Classical error models assume that an explanatory variable is measured with independent stochastic errors \cite{acom}. Berkson error models assume that the explanatory variable is a controlled variable with an error and that only the controlled variable can be measured \cite{one,burr}. Approaches to EIV models vary according to the situation. In this paper, we consider the former EIV. The corrected score function in Nakamura \cite{cor} has been used to estimate generalized linear models. In particular, the Poisson regression model is easy to handle analytically in generalized linear models as we see later. Thus, we focus on the Poisson regression model with measurement errors.

Approaches to a Poisson regression model with classical errors have been discussed by Kukush et al. \cite{three}, Shklyar and Schneeweiss \cite{com}, Jiang and Ma \cite{poi}, Guo and Li \cite{jie}, and so on. Kukush et al. \cite{three} described the statistical properties of the na\"{i}ve estimator, corrected score estimator, and structural quasi score estimator of a Poisson regression model with normally distributed explanatory variable and measurement error. Shklyar and Schneeweiss \cite{com} assumed an explanatory variable and a measurement error with a multivariate normal distribution and compared the asymptotic covariance matrices of the corrected score estimator, simple structural estimator, and structural quasi score estimator of a Poisson regression model. Jiang and Ma \cite{poi} assumed a high-dimensional explanatory variable with a multivariate normal error and proposed a new estimator for a Poisson regression model by combining Lasso regression and the corrected score function. Guo and Li \cite{jie} assumed a Poisson regression model with classical errors and proposed an estimator that is a generalization of the corrected score function discussed in Nakamura \cite{cor} for generally distributed errors; they derived the asymptotic normality of the proposed estimator.

In this study, we generalize the na\"{i}ve estimator discussed in Kukush et al. \cite{three}. They reported the bias of the na\"{i}ve estimator; however, the explanatory variable is not always normally distributed as they assume. In practice, the assumption of a normal distribution is not realistic. Here, we assume that the explanatory variable and measurement error are not limited to normal distributions. However, the na\"{i}ve estimator does not always exist in every situation. Therefore, we clarify the requirements for the existence of the na\"{i}ve estimator and derive its asymptotic bias. The constant vector to which the na\"{i}ve estimator converges in probability does not coincide with the unknown parameter in the model. Therefore, we propose a consistent estimator of the unknown parameter using the na\"{i}ve estimator. It is obtained from a system of equations that represent the relationship between the unknown parameter and constant vector. As illustrative examples, we present explicit representations of the new estimator for a Gamma explanatory variable with a normal error or a Gamma error.

In Section \ref{S2}, we present the Poisson regression model with measurement errors and the definition of the na\"{i}ve estimator and show that the na\"{i}ve estimator has an asymptotic bias for the true parameter. In Section \ref{S3}, we consider the requirements for the existence of the na\"{i}ve estimator and derive its asymptotic bias and asymptotic mean squared error (MSE) assuming that the explanatory variable and measurement error are generally distributed. In addition, we introduce application examples of a Gamma explanatory variable with a normal error or a Gamma error. In Section \ref{S4}, we propose the corrected na\"{i}ve estimator as a consistent estimator of the true parameter under general distributions and give application examples for a Gamma explanatory variable with a normal error or a Gamma error. In Section \ref{S5}, we present simulation studies that compare the performance of the na\"{i}ve estimator and corrected na\"{i}ve estimator.

\section{Preliminary}\label{S2}
In this section, we state the statistical model considered in this paper and the definition of the na\"{i}ve estimator and show that the na\"{i}ve estimator has an asymptotic bias for the true parameter.
\subsection{Poisson regression models with errors}
We assume a single covariate Poisson regression model between the objective variable $Y$ and explanatory variable $X$:
\begin{equation}
Y \mid X \sim Po(\exp(\beta_{0}+\beta_{1}X)). \nonumber 
\end{equation}
$X$ can typically be correctly observed. We assume here that $X$ has a stochastic error $U$ as
\begin{equation}
W = X + U, \nonumber 
\end{equation}
where $U$ is supposed to be independent of $(X,Y)$. We also assume that $(Y_i,X_i,U_i)~(i=1,\ldots,n)$ are i.i.d. samples of the distributions of $(Y \mid X,X,U)$. Although we can observe $Y$ and $W$, we assume that $X$ and $U$ cannot be directly observed. However, even if we know the family of the distributions of $X$ and $U$, we can-not make a statistical inference regarding $X$ and $U$ if we can observe only $W$. Because $U$ is the error distribution, the mean of $U$ is often zero, and we may suppose that we have empirical information about the degree of error (the variance of $U$). Therefore, in this study, we assume that the mean and variance of $U$ are known. From the above assumption, $Y$ and $W$ are independent for the given $X$.
\begin{align}
f_{Y,W \mid X}(y,w \mid x) &= \frac{f_{Y,W,X}(y,w,x)}{f_{X}(x)} = \frac{f_{Y,W,U}(y,w,w-x)}{f_{X}(x)} \nonumber \\
&= \frac{f_{Y,X}(y,x)f_{U}(w-x)}{f_{X}(x)} = f_{Y \mid X}(y \mid x)f_{W \mid X}(w \mid x). \nonumber
\end{align}
We use this conditional independence when we calculate the expectations.
\subsection{The \texorpdfstring{na\"{i}ve}{naive} estimator}
The na\"{i}ve estimator $\hat{\boldsymbol{\beta}}^{(N)}=(\hat{\beta}_{0}^{(N)},\hat{\beta}_{1}^{(N)})'$ for $\boldsymbol{\beta}=(\beta_0,\beta_1)'$ is defined as the solution of the equation
\begin{equation}
\label{1}
S_{n}(\hat{\boldsymbol{\beta}}^{(N)} \mid \mathcal{X}) = \boldsymbol{0}_{2},
\end{equation}
where
\begin{equation}
S_{n}(\boldsymbol{b} \mid \mathcal{X}) = \frac{1}{n} \sum_{i=1}^{n}\{Y_{i}- \exp(b_{0}+b_{1}W_{i})\}(1,W_i)', \nonumber 
\end{equation}
where $\boldsymbol{b} = (b_0,b_1)',\mathcal{X} = (X_1,\ldots,X_n)'$.
The na\"{i}ve estimator can be interpreted as the maximum likelihood estimator if we wrongly assume that $Y \mid W \sim Po(\exp(\beta_{0}+\beta_{1}W))$ because (\ref{1}) is the log-likelihood equation for $Y \mid W \sim Po(\exp(\beta_{0}+\beta_{1}W))$. The correct distribution of $Y \mid W$ is
\begin{align*}
f_{Y \mid W}(y \mid w) &= \frac{1}{f_{W}(w)}\int_{supp(f_{U})}f_{Y \mid W,U}(y \mid w,u)f_{U}(u)f_{X}(w-u) du  \\
&= \frac{1}{f_{W}(w)}\int_{supp(f_{U})}f_{Y \mid X}(y \mid w-u)f_{U}(u)f_{X}(w-u) du  \\
&= \frac{1}{f_{W}(w)}\int_{supp(f_{U})}Po(\exp(\beta_{0}+\beta_{1}(w-u)))f_{U}(u)f_{X}(w-u) du 
\end{align*}
assuming that $U$ is independent of $(X,Y)$. The right-hand side must be different from $Po(\exp(\beta_{0}+\beta_{1}W))$ in general. If one ignores the error $U$ and fits the likelihood estimation using $W$ instead of $X$, a biased estimator is obtained. In fact, by the law of large numbers, we have
\begin{align}
S_{n}(\hat{\boldsymbol{\beta}}^{(N)} \mid \mathcal{X}) &= \frac{1}{n} \sum_{i=1}^{n}\{Y_{i}- \exp(\hat{\beta}_{0}^{(N)}+\hat{\beta}_{1}^{(N)}W_{i})\}(1,W_i)' \nonumber \\
&\overset{p}{\longrightarrow}\textbf{E}_{X,Y \mid X,W}[\{Y - \exp(\hat{\beta}_{0}^{(N)}+\hat{\beta}_{1}^{(N)}W)\}(1,W)']. \nonumber 
\end{align}
Thus, the na\"{i}ve estimator converges to $\boldsymbol{b}=(b_0,b_1)'$ which is the solution of the estimating equation
\begin{equation}\label{alg1}
\textbf{E}_{X,Y \mid X,W}[\{Y - \exp(b_{0}+ b_{1}W)\}(1,W)'] = \boldsymbol{0}_{2}. 
\end{equation}
Eq. (\ref{alg1}) implies that for a given $\mathcal{X}$
\begin{equation}
\hat{\boldsymbol{\beta}}^{(N)} \overset{p}{\longrightarrow} \boldsymbol{b} \not= \boldsymbol{\beta}. \nonumber 
\end{equation}

\section{Properties of the \texorpdfstring{na\"{i}ve}{naive} estimator}\label{S3}
In this section, we consider the requirements for the existence of the na\"{i}ve estimator and derive its asymptotic bias and asymptotic MSE assuming that the explanatory variable and measurement error are generally distributed. In addition, we introduce application examples for a Gamma explanatory variable with a normal error or a Gamma error.

\subsection{The existence of the \texorpdfstring{na\"{i}ve}{naive} estimator}
The na\"{i}ve estimator does not always exist for general random variables $X$ and $U$. Thus, we assume the existence of the expectation
\begin{equation}
\textbf{E}_{X,Y \mid X,W}[\{Y - \exp(b_{0}+ b_{1}W)\}(1,W)'] \nonumber 
\end{equation}
as a requirement for the existence of the na\"{i}ve estimator. Consequently, the following four expectations should exist.
\begin{equation}
\label{2}
\begin{cases}
\textbf{E}[Y] &= \textbf{E}_{X}[\textbf{E}[Y \mid X]] = \textbf{E}_{X}[\exp(\beta_{0}+\beta_{1}X)] = e^{\beta_0}M_{X}(\beta_1) , \\
\textbf{E}[\exp(b_{0}+b_{1}W)] &= e^{b_0}\textbf{E}[e^{b_{1}X+b_{1}U}] = e^{b_0}M_{X}(b_1)M_{U}(b_1), \\
\textbf{E}[YW] &= \textbf{E}_{X}[\textbf{E}[Y \mid X]\textbf{E}[W \mid X]] = \textbf{E}_{X}[(X+\textbf{E}[U])\exp(\beta_{0}+\beta_{1}X)] \\
&= e^{\beta_0}\textbf{E}[U]M_{X}(\beta_1) + e^{\beta_0}\textbf{E}[Xe^{\beta_{1}X}]  \\
&= e^{\beta_0}\textbf{E}[U]M_{X}(\beta_1) + e^{\beta_0}\nabla M_{X}(\beta_{1}) , \\
\textbf{E}[W\exp(b_{0}+b_{1}W)] &= \textbf{E}_{X}[\textbf{E}_{U}[(X+U)\exp(b_{0}+b_{1}X+b_{1}U)]] \\
&= e^{b_0}\textbf{E}[Xe^{b_{1}X}]M_{U}(b_1) + e^{b_0}\textbf{E}[Ue^{b_{1}U}]M_{X}(b_1) \\
&= e^{b_0}M_{U}(b_1)\nabla M_{X}(b_{1}) + e^{b_0}M_{X}(b_1)\nabla M_{U}(b_{1}) .
\end{cases}
\end{equation}
Therefore, these expectations require that $M_{X}(\beta_1),M_{X}(b_1),M_{U}(b_1)$ exist. This condition is the requirement for the existence of the na\"{i}ve estimator. Here, we assume the existence of $M_{X}(\beta_1),M_{X}(b_1),M_{U}(b_1)$ for the distributions of $X$ and $U$.

\subsection{Asymptotic bias of the \texorpdfstring{na\"{i}ve}{naive} estimator}
The na\"{i}ve estimator satisfies
\begin{equation}
\hat{\boldsymbol{\beta}}^{(N)} \qquad  \overset{p}{\longrightarrow} \qquad \boldsymbol{b} \nonumber 
\end{equation}
and has an asymptotic bias for the true $\boldsymbol{\beta}$. Here, we derive the asymptotic bias under general conditions. From (\ref{alg1}), we obtain two equations:
\begin{equation}
\label{3}
\begin{cases}
\textbf{E}[Y] &= \textbf{E}[\exp(b_{0}+b_{1}W)], \\
\textbf{E}[YW] &= \textbf{E}[W\exp(b_{0}+b_{1}W)].
\end{cases}
\end{equation}
From (\ref{2}) with the above equalities, we have
\begin{align}
e^{\beta_0}M_{X}(\beta_1) &= e^{b_0}M_{X}(b_1)M_{U}(b_1), \nonumber \\
e^{\beta_0}\textbf{E}[U]M_{X}(\beta_1) + e^{\beta_0}\nabla M_{X}(\beta_1) &= e^{b_0}(\nabla M_{X}(b_1))M_{U}(b_1) + e^{b_0}(\nabla M_{U}(b_1))M_{X}(b_1) \nonumber \\
&= e^{b_0}\nabla (M_{X}(b_1)M_{U}(b_1)) =  e^{b_0}\nabla M_{W}(b_1). \nonumber
\end{align}
Therefore, we use a transformation to obtain the following system of equations:
\begin{align}
b_{0} &= \beta_{0}+\log \left(\frac{M_{X}(\beta_1)}{M_{W}(b_1)}\right) ,  \nonumber \\
K_{W}'(b_1) &=\frac{1}{M_{W}(b_1)}\nabla M_{W}(b_1) = \textbf{E}[U] + \frac{\nabla M_{X}(\beta_1)}{M_{X}(\beta_1)}, \nonumber 
\end{align}
where $K_{W}$ is the cumulant generating function of $W$. Thus, $\boldsymbol{b}=(b_{0},b_{1})'$ is determined by the solution of this system of equations. Therefore, the equation
\begin{equation}
K_{W}'(b_1)= \textbf{E}[U] + \frac{\nabla M_{X}(\beta_1)}{M_{X}(\beta_1)} \nonumber 
\end{equation}
should have a solution with respect to $b_1$. Here, we set
\begin{equation}
G(\beta_1,b_1):= K_{W}'(b_1)-\textbf{E}[U]-K_{X}'(\beta_1). \nonumber 
\end{equation}
We assume $G(\beta_1,b_1)$ has zero in $\mathbb{R}^{2}$ and satisfies
\begin{equation}
\frac{\partial G(\beta_1,b_1)}{\partial b_{1}}=K_{W}''(b_1) \not =0. \nonumber 
\end{equation}
Then, by the theorem of implicit functions, there exists a unique $C^{1}$-class function $g$ that satisfies $b_{1}=g(\beta_1)$ in the neighborhood of the zero of $G$. Using this expression, we write the asymptotic bias of the na\"{i}ve estimator as
\begin{align}
\lim_{n \to \infty}\textbf{E}[\hat{\beta}_{0}^{(N)}- \beta_{0}] &= b_{0}- \beta_{0} = \log\left(\frac{M_{X}(\beta_{1})}{M_{W}\circ g(\beta_{1})}\right), \nonumber \\
\lim_{n \to \infty}\textbf{E}[\hat{\beta}_{1}^{(N)}- \beta_{1}] &= b_{1}- \beta_{1} = g(\beta_{1})-\beta_{1}. \nonumber 
\end{align}
We also derive the asymptotic MSE of the na\"{i}ve estimator. The MSE can be represented as the sum of the squared bias and variance. The asymptotic variance of the na\"{i}ve estimator is $0$ because the na\"{i}ve estimator is a consistent estimator of $\boldsymbol{b}$. Thus, we obtain the asymptotic MSE of the na\"{i}ve estimator as
\begin{align}
\lim_{n \to \infty}\textbf{E}[(\hat{\beta}_{0}^{(N)}- \beta_{0})^{2}] &= (b_{0}- \beta_{0})^{2} = \left(\log\left(\frac{M_{X}(\beta_{1})}{M_{W}\circ g(\beta_{1})}\right)\right)^{2}, \nonumber \\
\lim_{n \to \infty}\textbf{E}[(\hat{\beta}_{1}^{(N)}- \beta_{1})^{2}] &= (b_{1}- \beta_{1})^{2} = (g(\beta_{1})-\beta_{1})^{2}. \nonumber 
\end{align}
Therefore, the asymptotic bias is given by the following theorem assuming general distributions.
\begin{theorem}[]\label{bias}
Let $Y \mid X \sim Po(\exp(\beta_{0}+\beta_{1}X))$. Assume that $W=X+U$ and $U$ is independent of $(X,Y)$. Assume the existence of $M_{X}(\beta_1),M_{X}(b_1),M_{U}(b_1)$. Let
\begin{equation}
G(\beta_1,b_1):= K_{W}'(b_1)-\textbf{E}[U]-K_{X}'(\beta_1). \nonumber 
\end{equation}
Assume G has zero in $\mathbb{R}^{2}$ and satisfies
\begin{equation}
\frac{\partial G(\beta_1,b_1)}{\partial b_{1}}=K_{W}''(b_1) \not =0. \nonumber 
\end{equation}
Then, the asymptotic biases of the na\"{i}ve estimators $\hat{\beta}^{(N)}_{0}$ and $\hat{\beta}^{(N)}_{1}$ are given by
\begin{equation}
\log\left(\frac{M_{X}(\beta_{1})}{M_{W}\circ g(\beta_{1})}\right), \quad g(\beta_{1})-\beta_{1} \nonumber 
\end{equation}
respectively, where $g$ is a $C^{1}$-class function satisfying $b_{1}=g(\beta_1)$ in the neighborhood of the zero of $G$. Furthermore, the asymptotic MSEs of the na\"{i}ve estimators $\hat{\beta}^{(N)}_{0}$ and $\hat{\beta}^{(N)}_{1}$ are given by their squared asymptotic biases.
\end{theorem}

\subsection{Examples}
In this section, we present two type of examples. First, we assume that a Gamma explanatory variable with a normal error. Let
\begin{equation}
X \sim \Gamma(k,\lambda),\quad U \sim N(0,\sigma^2), \nonumber 
\end{equation}
where $k>0, \lambda>0, 0<\sigma^2<\infty$. We apply the na\"{i}ve estimation under this condition. From the assumptions of Theorem \ref{bias}, we assume the existence of
\begin{equation}
M_{X}(\beta_1),M_{X}(b_1),M_{U}(b_1); \nonumber 
\end{equation}
therefore, we obtain the parameter conditions
\begin{equation}
\lambda - \beta_{1}>0, \quad \lambda - b_{1}>0. \nonumber 
\end{equation}
Next, we derive $\boldsymbol{b}=(b_0,b_1)'$. Under this condition, we obtain
\begin{equation}
G(\beta_{1},b_{1}) = K_{W}'(b_1)-\textbf{E}[U]-K_{X}'(\beta_1) = \frac{k}{\lambda-b_{1}}+\sigma^{2}b_{1}-\frac{k}{\lambda-\beta_{1}}. \nonumber 
\end{equation}
Thus, the set of zero of $G$ is
\begin{equation}
\left\{(\beta_{1},b_{1}) \in \mathbb{R}^{2} ; \beta_{1} = \frac{k+\lambda\sigma^{2}(\lambda-b_{1})}{k+\sigma^{2}(\lambda-b_{1})b_{1}}b_{1}\right\}. \nonumber 
\end{equation}
In addition,
\begin{equation}
\frac{\partial G(\beta_{1},b_{1})}{\partial b_{1}} = \frac{k}{(\lambda-b_{1})^{2}} + \sigma^{2} > 0. \nonumber 
\end{equation}
Therefore, $G$ has zero in $\mathbb{R}^{2}$ and satisfies $\frac{\partial G(\beta_{1},b_{1})}{\partial \beta_{1}} \not= 0$. From $G(\beta_{1},b_{1})=0$, we obtain two implicit functions
\begin{align}
b_{1}^{(1)} &= \frac{(\lambda-\beta_1)\lambda \sigma^{2}+k +\sqrt{s}}{2(\lambda-\beta_1)\sigma^2}, \nonumber \\
b_{1}^{(2)} &= \frac{(\lambda-\beta_1)\lambda \sigma^{2}+k -\sqrt{s}}{2(\lambda-\beta_1)\sigma^2}, \nonumber 
\end{align}
where $s = (\lambda-\beta_1)^{2}\lambda^{2} \sigma^{4} + 2(\lambda-\beta_1)(\lambda-2\beta_1)\sigma^{2}k + k^{2}>0$. Then, we obtain two expressions of $b_{0}$ corresponding to $b_{1}$.
\begin{align}
b_{0}^{(1)} &:= \beta_{0} + \log\left(\frac{M_{X}(\beta_1)}{M_{W}(b_{1}^{(1)})}\right) \nonumber \\
&= \beta_{0}  + k \log \frac{(\lambda- \beta_1)\lambda \sigma^{2} - k - \sqrt{s}}{2(\lambda- \beta_1)^{2}\sigma^2}   \nonumber \\
&\hspace{5mm}- \frac{(\lambda-\beta_1)^{2}\lambda^{2} \sigma^{4} + 2(\lambda-\beta_1)^{2}\sigma^{2}k + k^{2} +((\lambda-\beta_1)\lambda \sigma^{2} +k)\sqrt{s}}{4(\lambda-\beta_1)^{2}\sigma^{2}}, \nonumber \\
b_{0}^{(2)} &:= \beta_{0} + \log\left(\frac{M_{X}(\beta_1)}{M_{W}(b_{1}^{(2)})}\right) \nonumber \\
&= \beta_{0}  + k \log \frac{(\lambda- \beta_1)\lambda \sigma^{2} - k + \sqrt{s}}{2(\lambda- \beta_1)^{2}\sigma^2}   \nonumber \\
&\hspace{5mm}- \frac{(\lambda-\beta_1)^{2}\lambda^{2} \sigma^{4} + 2(\lambda-\beta_1)^{2}\sigma^{2}k + k^{2} -((\lambda-\beta_1)\lambda \sigma^{2} +k)\sqrt{s}}{4(\lambda-\beta_1)^{2}\sigma^{2}}.\nonumber 
\end{align}
In addition,
\begin{equation}
s = ((\lambda-\beta_{1})\lambda \sigma^{2}- k)^{2} + 4(\lambda-\beta_{1})^{2}\sigma^{2}k; \nonumber 
\end{equation}
therefore, $s$ satisfies $\sqrt{s}> \mid (\lambda-\beta_{1})\lambda \sigma^{2}- k \mid$. From the antilogarithm condition, $\boldsymbol{b} = (b_{0}^{(2)},b_{1}^{(2)})'$ is a solution of the system of equations (\ref{3}) in the range of $\mathbb{R}^{2}$. Thus, the asymptotic bias is given by
\begin{align}
b_{0}-\beta_{0} &= k \log \frac{(\lambda- \beta_1)\lambda \sigma^{2} - k + \sqrt{s}}{2(\lambda- \beta_1)^{2}\sigma^2}  \nonumber \\
&\hspace{5mm}- \frac{(\lambda-\beta_1)^{2}\lambda^{2} \sigma^{4} + 2(\lambda-\beta_1)^{2}\sigma^{2}k + k^{2} -((\lambda-\beta_1)\lambda \sigma^{2} +k)\sqrt{s}}{4(\lambda-\beta_1)^{2}\sigma^{2}}, \nonumber \\
b_{1}-\beta_{1} &= \frac{\lambda}{2}-\beta_{1} + \frac{k-\sqrt{s}}{2(\lambda-\beta_{1})\sigma^{2}}. \nonumber 
\end{align}

Next, we present another example, Gamma explanatory variable with a Gamma error. Let
\begin{equation}
X \sim \Gamma(k_1,\lambda),\quad U \sim \Gamma(k_2,\lambda), \nonumber 
\end{equation}
where $k_{1}>0,k_{2}>0,\lambda>0$. We apply the na\"{i}ve estimation under this condition. From the assumptions of Theorem \ref{bias}, we assume the existence of
\begin{equation}
M_{X}(\beta_1),M_{X}(b_1),M_{U}(b_1); \nonumber 
\end{equation}
therefore, we obtain the parameter conditions
\begin{equation}
\lambda - \beta_{1}>0, \quad \lambda - b_{1}>0 .\nonumber 
\end{equation}
Next, we derive $\boldsymbol{b}=(b_0,b_1)'$. Under this condition, we obtain
\begin{equation}
G(\beta_{1},b_{1}) = \frac{k_{1}+k_{2}}{\lambda-b_{1}}-\frac{k_{1}}{\lambda-\beta_{1}}-\frac{k_{2}}{\lambda}.\nonumber 
\end{equation}
Thus, the set of zero of $G$ is
\begin{equation}
\left\{(\beta_{1},b_{1}) \in \mathbb{R}^{2}; b_{1} = \frac{k_{1}\lambda\beta_{1}}{k_{1}\lambda+k_{2}(\lambda-\beta_{1})}\right\}. \nonumber 
\end{equation}
In addition,
\begin{equation}
\frac{\partial G(\beta_{1},b_{1})}{\partial b_{1}} = \frac{k_{1}+k_{2}}{(\lambda-b_{1})^{2}}>0. \nonumber 
\end{equation}
Therefore, $G$ has zero in $\mathbb{R}^{2}$ and satisfies $\frac{\partial G(\beta_{1},b_{1})}{\partial \beta_{1}} \not= 0$. From $G(\beta_{1},b_{1}) = 0$, we obtain the implicit function
\begin{equation}
b_{1} = \frac{k_{1}\lambda\beta_{1}}{k_{1}\lambda+k_{2}(\lambda-\beta_{1})}. \nonumber 
\end{equation}
Thus, by Theorem \ref{bias}, the asymptotic bias is given by
\begin{align}
b_{0}-\beta_{0} &= - k_{1}\log(1-\beta_{1}/\lambda) + (k_1 + k_2)\log(1-b_{1}/\lambda), \nonumber \\
b_{1}-\beta_{1} &= -\frac{k_{2}(\lambda-\beta_{1})\beta_{1}}{k_{1}\lambda+k_{2}(\lambda-\beta_{1})}. \nonumber 
\end{align}

\section{Corrected \texorpdfstring{na\"{i}ve}{naive} estimator}\label{S4}
In this section, we propose a corrected na\"{i}ve estimator as a consistent estimator of $\boldsymbol{\beta}$ under general distributions and give application examples for a Gamma explanatory variable with a normal error or a Gamma error. From (\ref{alg1}),
\begin{align}
\textbf{E}[Y] &= \textbf{E}[\exp(b_{0}+b_{1}W)], \nonumber \\
\textbf{E}[YW] &= \textbf{E}[W\exp(b_{0}+b_{1}W)], \nonumber 
\end{align}
we obtain the following system of equations:
\begin{align}
\beta_{0} &= b_{0} + \log \left(\frac{M_{W}(b_{1})}{M_{X}(\beta_{1})}\right),\nonumber \\
G(\beta_{1},b_{1}) &= K_{W}'(b_{1}) - \textbf{E}[U] - K_{X}'(\beta_{1}) = 0. \nonumber 
\end{align}
By solving this system of equations for $\beta_{0},\beta_{1}$ and replacing $\boldsymbol{b}=(b_{0},b_{1})'$ with the na\"{i}ve estimator $\hat{\boldsymbol{\beta}}^{(N)} = (\hat{\beta}_{0}^{(N)},\hat{\beta}_{1}^{(N)})'$, we obtain the consistent estimator of the true $\boldsymbol{\beta}$. Here,
\begin{equation}
\hat{\boldsymbol{\beta}}^{(N)} = \left(
	\begin{array}{c}
		\hat{\beta}_{0}^{(N)}  \\
		\hat{\beta}_{1}^{(N)}
		\end{array}
	\right) \qquad \overset{p}{\longrightarrow} \qquad \boldsymbol{b} = \left(
	\begin{array}{c}
		b_{0}  \\
		b_{1}
		\end{array}
	\right). \nonumber 
\end{equation}
Therefore,
\begin{equation}
\hat{\boldsymbol{\beta}}^{(CN)} \qquad \overset{p}{\longrightarrow} \qquad \boldsymbol{\beta}. \nonumber 
\end{equation}
Thus, $\hat{\boldsymbol{\beta}}^{(CN)}$ is a consistent estimator of $\boldsymbol{\beta}$. If $G$ has zero in $\mathbb{R}^{2}$ and satisfies
\begin{equation}
\frac{\partial G(\beta_1,b_1)}{\partial \beta_{1}} = -K_{X}''(\beta_1) \not= 0, \nonumber 
\end{equation}
then, by the theorem of implicit functions, there exists a unique $C^{1}$-class function $h$ that satisfies $\beta_{1}=h(b_1)$ in the neighborhood of the zero of $G$. We note that $h$ is the inverse function of $g$ in Theorem \ref{bias}. We propose a corrected na\"{i}ve estimator that is the consistent estimator of the true $\boldsymbol{\beta}$ as follows.
\begin{theorem}[]\label{CN}
Let $Y \mid X \sim Po(\exp(\beta_{0}+\beta_{1}X))$. Assume that $W=X+U$ and $U$ is independent of $(X,Y)$. Assume the existence of $M_{X}(\beta_1),M_{X}(b_1),M_{U}(b_1)$. Let
\begin{equation}
G(\beta_1,b_1) := K_{W}'(b_1)-\textbf{E}[U]-K_{X}'(\beta_1). \nonumber 
\end{equation}
Assume $G$ has zero in $\mathbb{R}^{2}$ and satisfies
\begin{equation}
\frac{\partial G(\beta_1,b_1)}{\partial \beta_{1}} = -K_{X}''(\beta_1) \not= 0. \nonumber 
\end{equation}
Then, the corrected na\"{i}ve estimator $\hat{\boldsymbol{\beta}}^{(CN)} = (\hat{\beta}_{0}^{(CN)},\hat{\beta}_{1}^{(CN)})'$, which corrects the bias of the na\"{i}ve estimator $\hat{\boldsymbol{\beta}}^{(N)} = (\hat{\beta}_{0}^{(N)},\hat{\beta}_{1}^{(N)})'$, is given by
\begin{align}
\hat{\beta}_{0}^{(CN)} &= \hat{\beta}_{0}^{(N)} + \log\left(\frac{M_{W}(\hat{\beta}_{1}^{(N)})}{M_{X}(\hat{\beta}_{1}^{(CN)})}\right), \nonumber \\
\hat{\beta}_{1}^{(CN)} &= h(\hat{\beta}_{1}^{(N)}), \nonumber 
\end{align}
where $h$ is a $C^{1}$-class function satisfying $\beta_{1}=h(b_1)$ in the neighborhood of the zero of $G$. Furthermore, the corrected na\"{i}ve estimator is a consistent estimator of $\boldsymbol{\beta}$.
\end{theorem}

\begin{example}
We derive the corrected na\"{i}ve estimator assuming
\begin{equation}
X \sim \Gamma(k,\lambda),U \sim N(0,\sigma^2). \nonumber 
\end{equation}
We obtain
\begin{align}
G(\beta_{1},b_{1}) &= \frac{k}{\lambda-b_{1}}+\sigma^{2}b_{1}-\frac{k}{\lambda-\beta_{1}}, \nonumber \\
\frac{\partial G(\beta_{1},b_{1})}{\partial \beta_{1}} &= -\frac{k}{(\lambda-\beta_{1})^{2}}<0. \nonumber 
\end{align}
$G$ has zero in $\mathbb{R}^{2}$ and satisfies $\frac{\partial G(\beta_{1},b_{1})}{\partial \beta_{1}} \not= 0$. From $G(\beta_{1},b_{1}) = 0$, we obtain the implicit function
\begin{equation}
\beta_{1} = \frac{\sigma^{2}\lambda b_{1}^{2} - (k+\lambda^{2}\sigma^{2})b_{1}}{\sigma^{2}b_{1}^{2}-\lambda\sigma^{2}b_{1}-k}= h(b_{1}). \nonumber 
\end{equation}
Thus, by Theorem \ref{CN}, the corrected na\"{i}ve estimator is given by
\begin{align}
\hat{\beta}_{0}^{(CN)} &= \hat{\beta}_{0}^{(N)} + \log\left(\frac{M_{W}(\hat{\beta}_{1}^{(N)})}{M_{X}(\hat{\beta}_{1}^{(CN)})}\right) \nonumber \\
&= \hat{\beta}_{0}^{(N)} + \frac{1}{2}\hat{\beta}_{1}^{(N)2}\sigma^{2} +k \log(1-\hat{\beta}_{1}^{(CN)}/\lambda)-k \log(1-\hat{\beta}_{1}^{(N)}/\lambda), \nonumber \\
\hat{\beta}_{1}^{(CN)} &= h(\hat{\beta}_{1}^{(N)}) = \frac{\lambda\sigma^{2}\hat{\beta}_{1}^{(N)2}-(k+\lambda^{2}\sigma^{2})\hat{\beta}_{1}^{(N)}}{\sigma^{2}\hat{\beta}_{1}^{(N)2}-\lambda\sigma^{2}\hat{\beta}_{1}^{(N)}-k}.\nonumber 
\end{align}
\end{example}

\begin{example}
We derive the corrected na\"{i}ve estimator assuming
\begin{equation}
X \sim \Gamma(k_{1},\lambda),U \sim \Gamma(k_{2},\lambda). \nonumber 
\end{equation}
We obtain
\begin{align}
G(\beta_{1},b_{1}) &= \frac{k_{1}+k_{2}}{\lambda-b_{1}}-\frac{k_{1}}{\lambda-\beta_{1}}-\frac{k_{2}}{\lambda}, \nonumber \\
\frac{\partial G(\beta_{1},b_{1})}{\partial \beta_{1}} &= -\frac{k_{1}}{(\lambda-\beta_{1})^{2}}<0. \nonumber 
\end{align}
$G$ has zero in $\mathbb{R}^{2}$ and satisfies $\frac{\partial G(\beta_{1},b_{1})}{\partial \beta_{1}} \not= 0$. From $G(\beta_{1},b_{1}) = 0$, we obtain the implicit function
\begin{equation}
\beta_{1} =\frac{(k_{1}+k_{2})b_{1}\lambda}{k_{1}\lambda+k_{2}b_{1}} = h(b_{1}). \nonumber 
\end{equation}
Thus, by Theorem \ref{CN}, the corrected na\"{i}ve estimator is given by
\begin{align}
\hat{\beta}_{0}^{(CN)} &= \hat{\beta}_{0}^{(N)} + \log\left(\frac{M_{W}(\hat{\beta}_{1}^{(N)})}{M_{X}(\hat{\beta}_{1}^{(CN)})}\right) \nonumber \\
&= \hat{\beta}_{0}^{(N)} + k_{1}\log(1-\hat{\beta}_{1}^{(CN)}/\lambda)-(k_{1}+k_{2})\log(1-\hat{\beta}_{1}^{(N)}/\lambda), \nonumber \\
\hat{\beta}_{1}^{(CN)} &= h(\hat{\beta}_{1}^{(N)}) = \frac{(k_{1}+k_{2})\hat{\beta}_{1}^{(N)}\lambda}{k_{1}\lambda+k_{2}\hat{\beta}_{1}^{(N)}}. \nonumber 
\end{align}
\end{example}

\section{Simulation studies}\label{S5}
In this section, we present simulation studies that compare the performance of the na\"{i}ve estimator and corrected na\"{i}ve estimator. We denote the sample size by $n$ and the number of simulations by MC. We calculate the estimated bias for $\hat{\boldsymbol{\beta}}^{(N)}$ and $\hat{\boldsymbol{\beta}}^{(CN)}$ as follows:
\begin{align}
\widehat{\mbox{BIAS}(\hat{\boldsymbol{\beta}}^{(N)})} &= \frac{1}{MC}\sum_{i=1}^{MC}\hat{\boldsymbol{\beta}}_{i}^{(N)} - \boldsymbol{\beta}, \nonumber \\
\widehat{\mbox{BIAS}(\hat{\boldsymbol{\beta}}^{(CN)})} &= \frac{1}{MC}\sum_{i=1}^{MC}\hat{\boldsymbol{\beta}}_{i}^{(CN)} - \boldsymbol{\beta},\nonumber
\end{align}
where $\hat{\boldsymbol{\beta}}_{i}^{(N)}$ and $\hat{\boldsymbol{\beta}}_{i}^{(CN)}$ represent the na\"{i}ve estimator and corrected na\"{i}ve estimator in the $i$th time simulation, respectively. Similarly, we calculate the estimated MSE matrix for $\hat{\boldsymbol{\beta}}^{(N)}$ and $\hat{\boldsymbol{\beta}}^{(CN)}$ as follows:
\begin{align}
\widehat{\mbox{MSE}(\hat{\boldsymbol{\beta}}^{(N)})} &= \frac{1}{MC}\sum_{i=1}^{MC}(\hat{\boldsymbol{\beta}}_{i}^{(N)}-\boldsymbol{\beta})(\hat{\boldsymbol{\beta}}_{i}^{(N)}-\boldsymbol{\beta})',\nonumber \\
\widehat{\mbox{MSE}(\hat{\boldsymbol{\beta}}^{(CN)})} &= \frac{1}{MC}\sum_{i=1}^{MC}(\hat{\boldsymbol{\beta}}_{i}^{(CN)}-\boldsymbol{\beta})(\hat{\boldsymbol{\beta}}_{i}^{(CN)}-\boldsymbol{\beta})'. \nonumber
\end{align}

\subsection{Case 1}
We assume $X \sim \Gamma(k,\lambda),U \sim N(0,\sigma^2)$. Let $\beta_{0}=0.2,\beta_{1}=0.3,k=2,\lambda=1.2,n=500,MC=1000$. We perform simulations with $\sigma^{2}=0.05,0.5,2$. Note that we assume that the true value of $\sigma^{2}$ is known. We estimate $k,\lambda$ in the formula of the corrected na\"{i}ve estimator by the moment method in terms of $W$ because the value of $X$ cannot be directly observed.
\begin{align}
\hat{k} &= \left(\frac{1}{n}\sum_{i=1}^{n}w_{i}\right) \hat{\lambda}, \nonumber\\
\hat{\lambda} &= \frac{\frac{1}{n}\sum_{i=1}^{n}w_{i}}{\frac{1}{n}\sum_{i=1}^{n}(w_{i}-\bar{w})^{2} - \sigma^{2}}.\nonumber
\end{align}

\begin{table}[h]
\begin{center}
\begin{minipage}{235pt}
\caption{Estimated bias of a Gamma distribution with a Normal error}\label{t1}%
\begin{tabular}{@{}llllll@{}}
\toprule
 &  & Asy.Bias $\hat{\beta}_{0}$ & $\widehat{\mbox{BIAS}(\hat{\beta}_{0})}$ & Asy.Bias $\hat{\beta}_{1}$ & $\widehat{\mbox{BIAS}(\hat{\beta}_{1})}$ \\
\midrule
\multirow{2}{*}{$\sigma^{2}=0.05$} & Na\"{i}ve & 0.01111 & 0.01139 & -0.005993 & -0.007199 \\
& CN & 0 & 0.00003532 & 0 & 0.0002603  \\
\multirow{2}{*}{$\sigma^{2}=0.5$} & Na\"{i}ve & 0.09912 & 0.1025 & -0.05297 & -0.05582  \\
& CN & 0 & 0.007817 & 0 & 0.0007142  \\
\multirow{2}{*}{$\sigma^{2}=2$} & Na\"{i}ve & 0.2757 & 0.2774 & -0.1454 & -0.1472  \\
& CN & 0 & -0.009493 & 0 & 0.002736  \\
\botrule
\end{tabular}
\end{minipage}
\end{center}
\end{table}
Table \ref{t1} shows the estimated bias of the true $\boldsymbol{\beta}$. Asy.Bias $\hat{\beta}_{0}$ and Asy.Bias $\hat{\beta}_{1}$ denote the theoretical asymptotic biases of $\hat{\beta}^{(N)}_{0}$ and $\hat{\beta}^{(N)}_{1}$, respectively, given in Theorem \ref{bias}. The bias correction of the na\"{i}ve estimator is done by the corrected na\"{i}ve estimator. With increasing $\sigma^{2}$, the bias of the na\"{i}ve estimator increases. However, the bias of the corrected na\"{i}ve estimator is small for large $\sigma^{2}$.
\begin{table}[h]
\begin{center}
\begin{minipage}{235pt}
\caption{Estimated MSE of a Gamma distribution with a Normal error}\label{t2}%
\begin{tabular}{@{}llllll@{}}
\toprule
 &  & Asy.MSE $\hat{\beta}_{0}$ & $\widehat{\mbox{MSE}(\hat{\beta}_{0})}$ & Asy.MSE $\hat{\beta}_{1}$ & $\widehat{\mbox{MSE}(\hat{\beta}_{1})}$ \\
\midrule
\multirow{2}{*}{$\sigma^{2}=0.05$} & Na\"{i}ve & 0.0001235 & 0.003003 & 0.00003592 & 0.0004536  \\
& CN & 0 & 0.002920 & 0 & 0.0004254  \\
\multirow{2}{*}{$\sigma^{2}=0.5$} & Na\"{i}ve & 0.009824 & 0.01362 & 0.002806 & 0.003508  \\
& CN & 0 & 0.003806 & 0 & 0.0006354  \\
\multirow{2}{*}{$\sigma^{2}=2$} & Na\"{i}ve & 0.07600 & 0.08124 & 0.02115 & 0.02214  \\
& CN & 0 & 0.01021 & 0 & 0.002160  \\
\botrule
\end{tabular}
\end{minipage}
\end{center}
\end{table}
Table \ref{t2} shows the estimated MSE of the true $\boldsymbol{\beta}$. Asy.MSE $\hat{\beta}_{0}$ and Asy.MSE $\hat{\beta}_{1}$ denote the theoretical asymptotic MSEs of $\hat{\beta}^{(N)}_{0}$ and $\hat{\beta}^{(N)}_{1}$, respectively, given in Theorem \ref{bias}. The MSE of the corrected na\"{i}ve estimator is smaller than that of the na\"{i}ve estimator in all cases.

\subsection{Case 2}
We assume $X \sim \Gamma(k_{1},\lambda),U \sim \Gamma(k_{2},\lambda)$. Let $\beta_{0}=0.2,\beta_{1}=0.3,k_{1}=2,\lambda=1.2,n=500,MC=1000$. We perform simulations with $k_{2}=0.072,0.72,2.88$. Similarly, we assume that the true value of $k_{2}$ is known. We estimate $k_{1},\lambda$ in the formula of the corrected na\"{i}ve estimator by the moment method in terms of $W$ because the value of $X$ cannot be directly observed.
\begin{align}
\hat{k}_{1} &= \left(\frac{1}{n}\sum_{i=1}^{n} w_{i}\right)\hat{\lambda}-k_{2},\nonumber \\
\hat{\lambda} &= \frac{\frac{1}{n}\sum_{i=1}^{n} w_{i}}{\frac{1}{n}\sum_{i=1}^{n}(w_{i}-\bar{w})^{2}}.\nonumber
\end{align}
\begin{table}[h t b]
\begin{center}
\begin{minipage}{235pt}
\caption{Estimated bias of a Gamma distribution with a Gamma error}\label{t3}%
\begin{tabular}{@{}llllll@{}}
\toprule
 &  & Asy.Bias $\hat{\beta}_{0}$ & $\widehat{\mbox{BIAS}(\hat{\beta}_{0})}$ & Asy.Bias $\hat{\beta}_{1}$ & $\widehat{\mbox{BIAS}(\hat{\beta}_{1})}$ \\
\midrule
\multirow{2}{*}{$k_{2}=0.072$} & Na\"{i}ve & -0.002634 & -0.005415 & -0.007887 & -0.008874  \\
& CN & 0 & -0.0006636 & 0 & 0.0002777  \\
\multirow{2}{*}{$k_{2}=0.72$} & Na\"{i}ve & -0.02090 & -0.01725 & -0.06378 & -0.06475  \\
& CN & 0 & -0.0002963 & 0 & -0.003184  \\
\multirow{2}{*}{$k_{2}=2.88$} & Na\"{i}ve & -0.04953 & -0.05439 & -0.1558 & -0.1569  \\
& CN & 0 & 0.002954 & 0 & -0.003224  \\
\botrule
\end{tabular}
\end{minipage}
\end{center}
\end{table}
Table \ref{t3} shows the estimated bias of the true $\boldsymbol{\beta}$. Similarly, the bias correction of the na\"{i}ve estimator is done by the corrected na\"{i}ve estimator. The bias of the corrected na\"{i}ve estimator is small when the variance of the error is large. Table \ref{t4} shows the estimated MSE of the true $\boldsymbol{\beta}$. The MSE of the corrected na\"{i}ve estimator is also smaller than that of the na\"{i}ve estimator.
\begin{table}[h]
\begin{center}
\begin{minipage}{235pt}
\caption{Estimated MSE of a Gamma distribution with a Gamma error}\label{t4}%
\begin{tabular}{@{}llllll@{}}
\toprule
 &  & Asy.MSE $\hat{\beta}_{0}$ & $\widehat{\mbox{MSE}(\hat{\beta}_{0})}$ & Asy.MSE $\hat{\beta}_{1}$ & $\widehat{\mbox{MSE}(\hat{\beta}_{1})}$ \\
\midrule
\multirow{2}{*}{$k_{2}=0.072$} & Na\"{i}ve & 0.08533 & 0.003109 & 0.000006940 & 0.0005384  \\
& CN & 0 & 0.003074 & 0 & 0.0004743 \\
\multirow{2}{*}{$k_{2}=0.72$} & Na\"{i}ve & 0.05580 & 0.005320 & 0.0004368 & 0.004894  \\
& CN & 0 & 0.004457 & 0 & 0.0008818  \\
\multirow{2}{*}{$k_{2}=2.88$} & Na\"{i}ve & 0.02080 & 0.01147 & 0.002453 & 0.02553  \\
& CN & 0 & 0.007401 & 0 & 0.001963  \\
\botrule
\end{tabular}
\end{minipage}
\end{center}
\end{table}

\section*{Declarations}
On behalf of all authors, the corresponding author states that there is no conflict of interest.


\bibliography{sn-bibliography}


\end{document}